\documentclass[60x84/16,10pt]{acmpt}
\usepackage[cp1251]{inputenc}

\usepackage[T1,T2A]{fontenc}
\usepackage{ucs}
\usepackage[english,russian]{babel}

\usepackage{amsmath}
\usepackage{amssymb}
\usepackage{amscd}

\usepackage{mathtools}
\mathtoolsset{
showonlyrefs=true,
mathic=true,
}

\allowdisplaybreaks

\usepackage{graphicx}
\usepackage{epstopdf}

\usepackage{hyperref}
\hypersetup{backref,
 colorlinks=false}
\hypersetup{pdfborder=0 0 0}

\usepackage{cite}

\begin{document}
\pagestyle{plain}

\title{Dimension of the space of invariant finitely additive measures of general Markov chains \\ and their ergodic properties\\}

\maketitle 

\begin{center}
 
{\bf A.I. Zhdanok}\protect\footnote[1]{{Department of Mathematical Analysis and MTM,\\
Tuvan State University, Kyzyl, Republic of Tuva, Russia},\\
{Laboratory of Mathematical Modeling,\\
Tuvan Institute for Exploration of Natural Resources of the Siberian Branch \\of the Russian Academy of Sciences, Kyzyl,  Republic of Tuva, Russia\\
	    e-mail: zhdanok@inbox.ru}}\\
	    
\end{center}

\date{}

{\bf Abstract:} General Markov chains with a countably additive tran\-sition probability in arbitrary phase space are considered. Markov operators extend from the space of countably additive measures to the space of finitely additive measures. In the author's papers a theorem was earlier proved that if all invariant finitely additive measures of a Markov chain are countably additive, i.e. there are no invariant purely finitely additive measures, then their subspace is finite-dimensional and the Markov chain satisfies the Doob-Doeblin quasicompactness conditions. In the same paper, a partial inversion of this theorem was proved with the dimension "one". In this paper we prove the inversion of this assertion for any finite dimensionality, but under certain additional conditions. The ergodic consequences are given. Examples and methods for studying their asymptotics with the aid of invariant purely finitely additive measures are given.\\



{\bf Keywords: }{ general Markov chains, an arbitrary phase space, Markov operators, invariant purely finitely additive measures, conditions for quasi\-compactness, ergodic properties.}


\newpage

\title{Размерность пространства инвариантных конечно-аддитивных мер общих цепей Маркова\\ и их эргодические свойства}

\author[1,2]{А. И. Жданок}

\address[1]{Кафедра математического анализа\\ и методики преподавания математики,\\
 Тувинский государственный университет,
  г.Кызыл, Россия}

\address[2]{Лаборатория математического моделирования,\\ 
Тувинский институт комплексного освоения природных ресурсов\\
СО РАН,  г.Кызыл, Россия\\
 e-mail: zhdanok@inbox.ru}

\begin{abstract}
    Рассматриваются общие цепи Маркова в произвольном фазовом пространстве. Марковские операторы продолжаются с  пространства счётно-аддитивных мер на пространство конечно-аддитивных мер. В работах автора ранее была доказана теорема о том, что, если все инвариантные конечно-аддитивные меры цепи Маркова счётно-аддитивны, т.е. отсутствуют инвариантные чисто конечно-аддитивные меры, то их подпространство конечномерно и цепь Маркова удовлетворяет условиям квазикомпактности Дуба-Деблина. Там же было доказано частичное обращение этой теоремы при размерности «единица». В настоящей работе доказывается обращение данного утверждения при любой конечной размерности, но при некоторых дополнительных условиях. Даны эргодические следствия. Приводятся примеры и методика изучения  их асимптотики с помощью инвариантных чисто конечно-аддитивных мер.\\
     
\end{abstract}

\keywords{общие цепи Маркова, произвольное фазовое пространство, Марковские операторы, инвариантные чисто конечно-адди\-тив\-ные меры, условия ква\-зиком\-пакт\-нос\-ти, эргодические свойства}

\maketitle 

\section{Введение}
\label{sec:intro}

Приведем используемые обозначения и сведения из [1] и [2].

Пусть $X$ – произвольное бесконечное множество и $\Sigma$ - сигма-алгебра его подмножеств, содержащая все одноточечные подмножества из $X$. Обозначим $B(X,\Sigma)$ - банахово пространство ограниченных $\Sigma$ -измеримых функций $f: X\to R$ с sup-нормой. Рассматриваются также банаховы пространства ограниченных мер $\mu : \Sigma\to R$ , с нормой, равной полной вариации меры $\mu$ (но можно использовать и эквивалентную sup-норму): $ba(X,\Sigma)$ - пространство конечно-аддитивных мер, $ca(X,\Sigma)$ - пространство счетно-аддитивных мер.  Если $\mu\ge{0}$, то $||\mu||=\mu(X)$.\medskip

Для топологических пространств $X$ используется также $C(X)$ - банахово пространство ограниченных непрерывных функций $f: X \to R$ с $sup$-нормой.

Конечно-аддитивная неотрицательная мера $\mu$ называется \underline{чис\-то} \underline{конеч\-но-ад\-ди\-тив\-ной} (чистым зарядом), если любая счетно-адди\-тив\-ная мера $\lambda$ , удовлетворяющая условию $0\le\lambda\le\mu$ , тождественно равна нулю. Любая конечно-аддитивная мера $\mu$ однозначно разлагается в сумму $\mu =\mu_1 + \mu_2$ , где $\mu_1$ - счетно-аддитивная, а $\mu_2$ - чисто конечно-аддитивная меры. Чисто конечно-аддитивные меры также образуют банахово пространство $pfa(X,\Sigma)$ с той же нормой, $ba(X,\Sigma) =ca(X,\Sigma) \oplus pfa(X,\Sigma)$. Если счетную аддитивность трактовать как непрерывность меры, как функции множеств, то чисто конечно-аддитивная мера $\mu$ - это ``разрывная мера''. Она существует тогда и только тогда, когда существует последовательность множеств $K_n\in\Sigma, n\in N$, $K_1\supset K_2\supset\dots$, $\lim K_n=\cap_{n=1}^{\infty} K_n=\O$, такая, что $\mu (K_n)\equiv\mu (X)$, т.е. $\lim\mu (K_n)=\mu(X)\ne 0=\mu(\lim K_n)=\mu(\O)$, [2].

Обозначим множества мер: 

$S_{ba}=\{\mu\in{ba(X,\Sigma)}:\mu\ge{0}, ||\mu||=1\},$
 
$S_{ca}=\{\mu\in{ca(X,\Sigma)}:\mu\ge{0}, ||\mu||=1\},$

$S_{pfa}=\{\mu\in{pfa(X,\Sigma)}:\mu\ge{0}, ||\mu||=1\}.$

Все меры из этих множеств будем называть \underline{вероятностными}.\medskip 

Цепи Маркова (ЦМ) на измеримом пространстве $(X,\Sigma)$ задаются своей переходной функцией (вероятностью) $p(x,E), x\in X, E\in\Sigma$, при обычных условиях: \medskip

1) $0\le p(x,E) \le{1}, p(x,X)=1$;

2) $p(\cdot,E)\in{B(X,\Sigma)}, \forall{E}\in\Sigma$;

3) $p(x,\cdot)\in{ca(X,\Sigma)}, \forall{x}\in{X}$. \medskip

Подчеркнем, что \underline{переходная функция у нас счетно-аддитивна} по второму аргументу, т.е. мы рассматриваем \underline{классические ЦМ.} \medskip

Переходная функция порождает два марковских линейных ограниченных положительных интегральных оператора:\medskip
 
$T: B(X,\Sigma)\to{B(X,\Sigma)}, (Tf)(x)=Tf(x)=\int_X f(y)p(x,dy),$

$\forall{f\in{B(X,\Sigma)}},\forall{x}\in{X};$ \medskip

$A: ca(X,\Sigma)\to{ca(X,\Sigma)}, (A\mu)(E)= A\mu(E)=\int_X p(x,E)\mu(dx),$

$\forall{\mu\in{ca(X,\Sigma)}},\forall{E}\in\Sigma.$ \medskip

Пусть начальная мера $\mu_0\in{S_{ca}}$. Тогда итерационная последовательность счетно-аддитивных вероятностных мер $\mu_{n}=A\mu_{n-1}\in{S_{ca}},n\in{N}$, обычно и отождествляется с цепью Маркова.\medskip

Принято считать, что подобный операторный подход в теории цепей Маркова впервые детально был разработан в работе Yosida и Kakutani [3] (1941 год), что и отражено в ее названии. Важно подчеркнуть, что в [3] все используемые меры, в том числе переходная функция $p(x,\cdot)$, счетно-аддитивны. Описанное ниже продолжение оператора A на пространство конечно-аддитивных мер было произведено рядом авторов намного позже.

Интересно также отметить, что один из соавторов работы [3] Yosida, позже (1952 год) стал одним из соавторов работы [2] (совместно с Hewitt), в которой разработаны основы современной общей теории конечно-аддитивных мер (конечно, на базе уже существующих работ других авторов). Однако, цепи Маркова, как возможный объект приложения новой теории, в [2] даже и не упонимаются, хотя Yosida уже являлся общепризнанным специалистом и в этой области. 

В настоящей статье, как и в предыдущих работах автора [4], [5] и других, изучаются некоторые проблемы, лежащие в области соприкосновения операторной теории общих цепей Маркова и теории конечно-аддитивных мер. \medskip

 Топологически сопряженным к пространству $B(X,\Sigma)$ является (изоморфно) пространство конечно-аддитивных мер: $B^*(X,\Sigma)=ba(X,\Sigma)$. При этом топологически сопряженным к оператору Т служит оператор $T^*:ba(X,\Sigma)\to{ba(X,\Sigma)}$, однозначно определяемый по известному правилу через интегральные ``скалярные произведения'':          

$$
\langle T^*\mu, f \rangle=\langle \mu, Tf\rangle  \text{ для всех } f\in{B(X,\Sigma)} \text{ и } \mu\in{ba(X,\Sigma)}.
$$ 

Оператор $T^*$  является продолжением оператора А на все пространство $ba(X,\Sigma)$ с сохранением его аналитического вида, т.е.  
$$
T^*\mu(E)=\int_X p(x,E)\mu(dx), \forall{\mu\in{ba(X,\Sigma)}}, \forall E\in\Sigma.
$$ 

Интеграл от ограниченной измеримой функции по ограниченной конечно-аддитивной мере строится по той же схеме, как и интеграл Лебега по произвольной ограниченной счетно-аддитивной мере (у нас как раз и функции и меры ограниченные). Оператор $T^*$ имеет собственное инвариантное подпространство $ca(X,\Sigma)$, т.е. $T^*[ca(X,\Sigma)] \subset{ca(X,\Sigma)}$, на котором он совпадает с изначальным оператором $A$. Теперь конструкция операторов Т и Т* уже функционально замкнута. Мы будем по-прежнему обозначать оператор Т* как А. 

В такой постановке естественно допустить к рассмотрению и марковские последовательности вероятностных конечно-аддитивных мер: 

$$
\mu_0\in{S_{ba}}, \mu_n = A\mu_{n-1}\in{S_{ba}}, n\in{N},
$$ 
сохраняя счетную аддитивность переходной функции $p(x,\cdot)$ по второму аргументу. Несмотря на это обстоятельство, образ $A\mu$ чисто конечно-аддитивной меры $\mu$  может остаться чисто конечно-аддитивным, т. е., вообще говоря, 
$A[ba(X,\Sigma)] \not \subset{ca(X,\Sigma)}$. \medskip

Допустимо и кардинально изменить постановку задачи - разрешить самой переходной функции $p(x,\cdot)$ быть всего лишь конечно-аддитивной мерой. Такие ЦМ тоже изучаются (в т.ч. в нашей работе [5]) и называются \underline{``конечно-аддитивными ЦМ''}, но в настоящей статье не рассматриваются. Таким образом, в нашем случае уместна следующая терминология: изучаются \underline{счетно-аддитивные} цепи Маркова, заданные на пространстве \underline{конечно-аддитивных мер}.\medskip

Обозначим множества инвариантных вероятностных мер ЦМ:

$\Delta_{ba}=\{\mu\in{S_{ba}}: \mu=A\mu\}$, 

$\Delta_{ca}=\{\mu\in{S_{ca}}: \mu=A\mu\}$, 

$\Delta_{pfa}=\{\mu\in{S_{pfa}}: \mu=A\mu\}$. \\

Классическая (счетно-аддитивная) цепь Маркова может иметь инвариантные вероятностные счетно-аддитивные меры, а может и не иметь, т.е. возможно $\Delta_{ca}=\O$ (например, у симметричного блуждания на Z). 

В работе [6] (1962 год) (Теорема 2.2) \v{S}idak доказал, что любая счетно-аддитивная ЦМ на произвольном измеримом пространстве $(X,\Sigma)$, продолженная на пространство конеч\-но-аддитивных мер, имеет хотя бы одну инвариантную конечно-аддитивную меру, т.е. всегда $\Delta_{ba}\ne\O$. Этот результат был затем доказан более простым способом в работе автора [4]. 

\v{S}idak в [6] (Теорема 2.5) также установил в общем случае, что, если конечно-аддитивная мера $\mu$ инвариантна $A\mu=\mu$, и $\mu=\mu_1+\mu_2$ - ее разложение на счетно-аддитивную и чисто конечно-аддитивную компоненты, то каждая из них также инвариантна: $A\mu_1=\mu_1$, $A\mu_2=\mu_2$. Следовательно, $\Delta_{ba}=co\{\Delta_{ca},\Delta_{pfa}\}$, и достаточно изучать инвариантные меры из $\Delta_{ca}$ и из $\Delta_{pfa}$, по отдельности. Здесь в записи $co\{\Delta_{ca},\Delta_{pfa}\}$ подразумевается выпуклая линейная оболочка множеств $\Delta_{ca}$ и $\Delta_{pfa}$. 

Изучением инвариантных конечно-аддитивных мер для ЦМ, заданных на \underline{топологических пространствах} $(X,\Sigma)$, занимался также ряд авторов (прежде всего - Foguel, см., например, [7]). Некоторая библиография с комментариями по-этому вопросу дается в [4] и [5].\medskip

Если размерность $dim\Delta = n<\infty$ (индексы у $\Delta$ опускаем), то существует алгебраический нормированный конечный базис Гамеля $\{\mu_1,\mu_2,\dots,\mu_n\}$ из мер соответствующего типа. Тогда $\Delta=co\{\mu_1,\mu_2,\dots,$ $\mu_n\}$. В $\Delta_{ca}$ базис можно выбрать из попарно сингулярных мер. В $\Delta_{ba}$ базис можно выбрать из мер, всего лишь попарно дизъюнктных в структурном смысле (см. [2], [4]). 

Напомним, что меры $\mu_1$ и $\mu_2\in{S_{ba}}$ называются \underline{сингулярными}, если существуют множества $K_1,K_2\in\Sigma$ такие, что $K_1\cap K_2=\O$ и $\mu_1(K_1)=1$, $\mu_2(K_2)=1$. Дизъюнктность конечно-аддитивных мер мы здесь не используем. Заметим лишь, что для счетно-аддитивных мер их дизъюнктность совпадает с сингулярностью (см. [2], пункты 1.11, 1.21 и др.). \medskip

Пусть дано семейство мер любого типа $\{\mu_\alpha\}_{\alpha \in I}$, где $I -$ произвольное семейство индексов, и $\mu_\alpha \in {S_{ba}}, \alpha \in I$. Будем говорить, что это множество мер \underline{попарно сингулярно}, если для любых $\alpha, \beta \in I, \alpha \neq \beta$, существуют такие множества $K_\alpha, K_\beta \in \Sigma$, что $\mu_\alpha(K_\alpha)=1$, $\mu_\beta(K_\beta)=1$ и $K_\alpha \cap K_\beta =\emptyset$. Множества $K_\alpha$ и $K_\beta$  здесь не однозначны. Их обычно называют множествами полной меры. Можно грубо называть их также ``носителями'' мер $\mu_\alpha$ и $\mu_\beta$. Однако, для топологических пространств $(X, \mathfrak{B})$  есть точное и однозначное определение ``носителя'' счетно-аддитивной меры, которое мы здесь не используем.

Отметим, что множества полной меры в определении попарной сингулярности $K_\alpha$ и $K_\beta$ зависят от второй меры, т. е. можно уточнить обозначения:
$$K_\alpha=K_\alpha^\beta ,  K_\beta=K_\beta^\alpha .$$

Будем говорить, что семейство мер (любого типа) $\{\mu_\alpha\}_{\alpha \in I}$ \underline{сингуляр-} \underline{но в совокупности}, если для любых двух подсемейств мер $\{\mu_\alpha\}_{\alpha \in I_1}$ и $\{\mu_\beta\}_{\beta \in I_2}$, $I_1, I_2 \subset I$ и $I_1 \cap I_2 = \emptyset$, существуют два множества $K_1, K_2 \in \Sigma$, такие, что $K_1 \cap K_2 = \emptyset$ и $\mu_\alpha(K_1)=1$, при всех $\alpha \in I_1$ и $\mu_\beta(K_2)=1$ при всех $\beta \in I_2$.

Очевидно, из сингулярности в совокупности следует попарная сингулярность семейства мер. Однако, обратное, вообще говоря, неверно. Но для конечного семейства мер $\{ \mu_1, \mu_2, \dots , \mu_n\}$ эти понятия совпадают.

Мы будем говорить, что бесконечное семейство мер $M$ линейно независимо, если любое его конечное подмножество линейно независимо. Понятно, что, если множество мер $M$ состоит из попарно сингулярных мер, то оно линейно независимо.

\section{Теоремы об инвариантных мерах и квазикомпактности} 
\label{sec:base-section}

В исследованиях ряда авторов (см., в том числе, [4], [5]) было показано, что размерность и состав множества \underline{инвариантных} \underline{ конечно-адди\-тив\-ных мер} марковского оператора А тесным образом связаны с одним из центральных вопросов эргодической теории ЦМ, а именно, с известными условиями квазикомпактности  марковских операторов (см., например, [3] и [8], Глава V, \S 5).

Напомним, что оператор А называется  \underline{квазикомпактным}, если существует \underline{компактный} оператор $A_1$ (переводящий каждое ограниченное множество в предкомпактное) и целое число $k\ge1$ такие, что $||A^k-A_1||<1$. В этом случае и саму ЦМ будем называть квазикомпактной. Достаточным, а в некоторых случаях и необходимым, условием квазикомпактности оператора А является общеизвестное условие Дуба-Деблина (D):  \\ 
$
(D)
\begin{cases}
\text{{\it существуют мера $\varphi\in{ca(X,\Sigma)}$, $\varphi\ge{0},\varepsilon>0$ и $k\ge1$ такие, что}}\\
\text{{\it из $\varphi(E)\le\varepsilon, E\in\Sigma$, следует $p^k(x,E)\le{1-\varepsilon}$ для всех $x\in{X}$.}}
\end{cases}
$ \medskip

Замечание: верхний индекс $k$ в $p^k$ означает порядок интегральной свертки переходной функции, а не её степень.

Если это условие выполнено, т.е. если оператор А квазикомпактен, то Цепь Маркова имеет конечное число сингулярных инвариантных счетно-аддитивных мер и эргодические средние равномерно сходятся в некотором смысле к ним в метрической топологии.

В литературе условие (D) нередко формулируют и в других формах (не всегда эквивалентных), в том числе в книге одного из авторов условия (D) Дуба [8], в работе [3], а также в книге Невё [9] (Глава V, пункт 5.5), где утверждается, что условие (D) всегда необходимо и достаточно для квазикомпактности марковских операторов. В нашей работе [5] также дается близкое условие (\~D), являющееся некоторой модификацией условия (D). Приведем это условие, но вначале введем следующие обозначения.\medskip

Пусть на $(X,\Sigma)$ задана ЦМ с переходной функцией $p(x,E)$ и марковскими операторами Т и А. Для любого $m\in N$ определим новую ЦМ с переходной функцией $q_m(x,E)$ и марковскими операторами $T_m$ и $A_m$ по правилам построения средних по Чезаро: 
$$q_m(x,E)=\frac{1}{m}\sum\limits_{k=1}^{m}p^k(x,E),\qquad T_m=\frac{1}{m}\sum\limits_{k=1}^{m}T^k,\qquad A_m=\frac{1}{m}\sum\limits_{k=1}^{m}A^k.$$

Такую ЦМ назовем \underline{конечно-осредненной ЦМ} (по исходной ЦМ).

Сформулируем новое условие (\~D): \medskip \\
$
(\tilde D)
\begin{cases}
\text{{\it существуют $\varphi\in ca(X,\Sigma)$, $\varphi\ge 0$, $\varepsilon>0$ и $m\ge 1$ такие, что}}\\
\text{{\it из $\varphi(E)<\varepsilon$, $E\in\Sigma$, следует $q_m(x,E)\le{1-\varepsilon}$ для всех $x\in X$.}}
\end{cases}
$ \medskip

Очевидно, что (\~D) является условием Дуба-Деблина (D) для конечно-осредненнной ЦМ (при фиксированном $m\ge 1$) с параметром k=1. Следовательно, если выполнено условие (\~D), то операторы $T_m$ и $A_m$ являются квазикомпактными.

В работе [5] (Теорема 12.1) показывается, что если выполнено условие (D), то выполнено и условие (\~D). Это означает, что если исходная ЦМ квазикомпактна, то квазикомпактны и все её конечно-осредненные ЦМ, начиная с некоторого номера $m\ge 1$.

Обозначим $\tilde\Delta_m$ семейство всех нормированных положительных конечно-аддитивных инвариантных мер для конечно-осредненной ЦМ с параметром $m$. Очевидно, что $\Delta_{ba}\subset\tilde\Delta_m$ при всех $m\in N$, и, возможно, $\Delta_{ba}\ne\tilde\Delta_m$ при наличии циклических подклассов.\medskip

В работе автора [5] была доказана следующая теорема:\medskip

\underline{Теорема 12.2 [5]}. {\it Для произвольной ЦМ условие (\~D) эквивалентно условию (*):
$$(*) \qquad \Delta_{ba}\subset{ca(X,\Sigma)},$$
которое означает, что все инвариантные конечно-аддитивные меры исходной ЦМ являются счетно-аддитивными, или, другими словами, исходная ЦМ не имеет инвариантных чисто конечно-аддитивных мер.}\medskip

Привлекательность условия (*) в том, что в нем, в отличие от условия Дуба-Деблина (D), нет никакой аналитики и динамики, а есть только «качественное» понятное утверждение.\medskip

Из Теоремы 12.2 [5], свойств квазикомпактной ЦМ и включения $\Delta_{ba}\subset\tilde\Delta_{m}$ сразу вытекает следующее утверждение.\medskip

\underline{Теорема 8.2} [4]. {\it Для произвольной ЦМ, если выполнено условие (*), т.е. если $\Delta_{ba}\subset{ca(X,\Sigma)}$, то размерность $dim\Delta_{ba}=n<\infty$.}\medskip

Но интересным оказалось то, что утверждение Теоремы 8.2 можно доказать и без привлечения конечно-осредненных ЦМ и условия (\~D), т.е. получить его не как простое следствие Теоремы 12.2, а как самостоятельный общий факт. Это и было доказано автором в [4] (Теоремы 8.1 и 8.2) с использованием слабых топологий и техники Банаховых пределов.

Наконец, возникло интуитивное предположение, что утверждение Теоремы 8.2 можно обратить. В той же работе [4] было доказано такое обращение, но только для случая размерности $n=1$:\medskip

\underline{Теорема 8.3} [4] {\it Для произвольной ЦМ, если $dim\Delta_{ba}=1$, т.е., если ЦМ имеет в $S_{ba}$ единственную инвариантную конечно-адди\-тив\-ную меру, то $\Delta_{ba}\subset{ca(X,\Sigma)}$, т.е. эта мера счетно-аддитивна.}

Доказательство этого факта также приводится без использования конечно-осредненной ЦМ.

При этом, т.е при $n=1$, выполнено условие Дуба-Деблина (D) (а также условие (\~D)), ЦМ квазикомпактна, и у нее нет инвариантных чисто конечно-аддитивных мер. 

В настоящей работе мы приводим и доказываем обращение Теоремы 8.2 [4], т.е. обобщаем Теорему 8.3 [4] уже для произвольной размерности $n\in{N}$, но при дополнительном условии $(\alpha)$:\\
$ 
(\alpha)
\begin{cases}
\text{{\it Пусть $\mu\in\Delta_{ba}$, и дано множество $K_{\mu}\in\Sigma$, такое, что  $\mu(K_{\mu})=1$.}}\\
\text{{\it Тогда существует множество $K\subset{K_{\mu}}$, такое, что $\mu(K)=1$ и оно}}\\
\text{{\it стохастически замкнуто, т.е. $p(x,K)=1$ для любого $x\in{K}$.}}
\end{cases}
$ \medskip

Очевидно, что при  $K_{\mu}=X$ условие $(\alpha)$ тривиально выполняется при $K=K_{\mu}=X$.

Если заранее предположить, что мера $\mu\in\Delta_{ba}$ счетно-аддитивна, то условие $(\alpha)$ также будет выполнено [4] (Теорема 6.2). \medskip

{\bf\underline{Теорема 1}.} {\it Пусть $dim\Delta_{ba}=n<\infty$, $\Delta_{ba}=co\{\mu_1,\mu_2,\dots,\mu_n\}$. Пусть все базисные меры $\mu_i$ из  $\Delta_{ba}$ попарно сингулярны, и для каждой из них выполнено условие $(\alpha)$.

Тогда  $\Delta_{ba}\subset{ca(X,\Sigma)}$, т.е. все инвариантные конечно-аддитив\-ные меры цепи Маркова являются счетно-аддитивными.} \medskip

{\it\underline{Доказательство.}} Из попарной сингулярности конечного семейства  мер, совпадающей с сингулярностью в совокупности, следует, что существуют попарно непересекающиеся множества $K_1,K_2,$ $\dots,K_n\in\Sigma$ такие, что $\mu_i(K_i)=1$ для  $i=1,2,\dots,n$, и $K_i\cap{K_j}=\O$ при $i\ne{j}$.

По условию $(\alpha)$, для каждой меры $\mu_i$, $i=1,\dots,n$, существует множество $K^i\subset{K_i}$ такое, что $\mu_i(K^i)=1$ и $p(x,K^i)=1$ для всех $x\in{K^i}$.

Образуем подпространства $(X_i,\Sigma_i)$ пространства $(X,\Sigma)$ по правилу $X_i=K^i\cap{X}$, $\Sigma_i=K^i\cap\Sigma$ для $i=1,2,\dots,n$. Сузим все меры $\mu_i$ на подпространства  $(X_i,\Sigma_i)$.

Сузим переходную функцию $p(x,E)$ на каждое подпространство $(X_i,\Sigma_i)$ по тождественному правилу: $\forall{x}\in{K^i}$, $\forall{E}\in\Sigma_i$, $p_i(x,E)=p(x,E)$. Легко проверить, что все меры $\mu_i$ будут единственными инвариантными конечно-аддитивными мерами для частных цепей Маркова, порождаемыми переходными функциями $p_i(x,E)$ на $(X_i,\Sigma_i)$ при $i=1,2,\dots,n$.

Теперь можно считать, что на каждом измеримом пространстве $(X_i,\Sigma_i)$ задана своя цепь Маркова с единственной инвариантной конечно-ад\-ди\-тив\-ной мерой $\mu_i$. Следовательно, по Теореме 8.3 [4], все инвариантные меры $\mu_i$ счетно-аддитивны на измеримых пространствах $(X_i,\Sigma_i)$ соответственно.

Теперь продолжим эти меры $\mu_i$ до мер $\widetilde{\mu_i}$ на все исходное измеримое пространство $(X,\Sigma)$ по правилу $\widetilde{\mu_i}(E)=\mu_i(E\cap{K^i})$ для каждого $E\in\Sigma$, т.е. "нулями". Очевидно, счетно-аддитивность мер $\widetilde{\mu_i}$ при этом сохраняется и они будут инвариантны уже на всем $(X,\Sigma)$, те. $\widetilde{\mu_i}\in{\Delta_{ca}}$ при всех $i=1,2,\dots,n$, и $\Delta_{ba}=\Delta_{ca}=co\{\widetilde{\mu_1},\widetilde{\mu_2},\dots,\widetilde{\mu_n}\}$. Теорема доказана. \medskip

{\bf\underline{Теорема 2}.} {\it Пусть $dim\Delta_{ba}=\infty$. Тогда ЦМ имеет инвариантные чисто конечно-аддитивные меры, т.е. $\Delta_{pfa} \neq \emptyset$ и $dim \Delta_{pfa} = \infty$.} \medskip

{\it\underline{Доказательство.}} Рассмотрим два возможных варианта для  $\Delta_{ca} \subset \Delta_{ba}$.
Пусть $dim\Delta_{ca}=n<\infty$, включая случай $\Delta_{ca} = \emptyset$. Тогда, поскольку $\Delta_{ba}=co (\Delta_{ca}, \Delta_{pfa})$, то $\Delta_{pfa} \neq \emptyset$  и  $dim \Delta_{pfa} = \infty$, и теорема доказана. 
 
Пусть $dim \Delta_{ca} = \infty$. Тогда, согласно нашей Теореме 8.1. ([4], §8, стр. 81) $\Delta_{pfa} \neq \emptyset$  и  $dim \Delta_{pfa} = \infty$. Теорема доказана.

Теорема 8.1. доказывается в  [5]  в технике Банаховых пределов  и $\tau_{B}$-слабой топологии.

В Теореме 2 не используется условие ($\alpha$) для мер из $\Delta_{ba}$. Однако, использование условия ($\alpha$) в Теореме 1 не позволяет объединить обе теоремы и сформулировать их как еще одно необходимое и достаточное условие существования (или не существования) инвариантных чисто конечно-аддитивных мер.\medskip

{\bf\underline{Следствие.}} {\it Пусть выполнены условия Теоремы 1. Тогда для ЦМ выполнено условие Дуба-Деблина (D) и условие (\~D), наше условие (*), у неё нет инвариантных чисто конечно-аддитивных мер, ЦМ является квазикомпактной, и обладает всеми соответствующими эргодическими свойствами.} \medskip

В книге Revuz [10] (Глава 6, §3, Теорема 3.5 и 3.7, стр. 240-243 в русском издании) доказывается утверждение, как отмечает автор, взятое у Horowitz [11], также о связи чисто конечно-аддитивных инвариантных мер с квазикомпактностью ЦМ. В английском издании [10] конечно-аддитивные меры называются ``mean'', а чисто конечно-аддитивные - ``pure mean'', в русском издании [10] эти меры называются ``среднее'' и ``чистое среднее'' соответственно. Все рассмотрение в [10] проводится для сепарабельных измеримых пространств $(X,\Sigma)$, в которых сигма-алгебра $\Sigma$ порождена счетным семейством множеств (у нас сигма-алгебра $\Sigma$ произвольна). Заранее предполагается, что ЦМ удовлетворяет условию Харриса ([10], Глава 3, §2, Определение 2.6) с некоторой уже заданной инвариантной счетно-аддитивной $\sigma$-ограниченной мерой $\varphi$. У Horowitz [11] изначально предполагается, что переходная функция $p(x,\cdot)$ при каждом $x \in X$ абсолютно непрерывна относительно заранее заданной счётно-аддитивной меры $m$.

При этих жёстких предположениях доказывается, что условие –{\it ЦМ является квазикомпактной}, эквивалентно условию – {\it не существует инвариантной чисто конечно-аддитивной меры} для ЦМ ([10], Теоремы 3.5 и 3.7). Утверждается, что при выполнении этих эквивалентных условий заранее заданная инвариантная счетно-аддитивная мера $\varphi$ из условия Харриса оказывается ограниченной и \underline{единственной} инвариантной конечно-аддитивной мерой для ЦМ ([10], Теорема 3.7).\medskip

В работе Horowitz [11] (пункт 4, теорема 4.1) доказывается аналогичное утверждение, но для марковского оператора $P$, заданного слева и справа на пространствах $L_\infty (m)$ и ${L_\infty}^* (m)$, соответственно. Для конечно-аддитивных мер используются термины ``charge'' и ``pure charge'' (заряд).\medskip

В наших рассмотрениях ЦМ не предполагается харрисовской, и она вообще не привязана к какой-либо заранее заданной и уже  инвариантной мере $m$. Наши Теоремы 1,2,  а также Теоремы 8.2, 8.3 из [4] и Теорема 12.2 из [5] носят гораздо более общий характер и содержат более сильные утверждения, чем в работах [10] и [11].\medskip

Для решения подобных задач в работах автора [4], [5] был разработан соответствующий аппарат. Сейчас мы воспользуемся следующей Теоремой 7.8., доказанной в [4] (Глава 1, \S 7, стр. 80), дающей аналитические условия (критерий Z) существования инвариантных чисто конечно-аддитивных мер для общих ЦМ на произвольных фазовых пространствах. Приведем эти условия, поскольку мы ими будем пользоваться  в последующих Примерах. \medskip

\underline{Теорема 7.8. [4].} {\it Пусть на произвольном измеримом пространстве $(X,\Sigma)$ задана счётно-аддитивная ЦМ с переходной функцией $p(x,E)$, удовлетворяющей следующему условию:} \medskip

$
({\bf{Z}})
\begin{cases}
\text{{\it Существуют последовательности чисел $\varepsilon_n$ и множеств $K_n$,}}\\
\text{{\it такие, что $\varepsilon_n\ge 0$, $\varepsilon_n\to 0$ при $n\to\infty$, $K_n\in\Sigma$, $K_n\ne\O$}}\\
\text{{\it при $n\in N$, $K_1\supset K_2\supset\dots, \bigcap\limits_{n=1}^\infty K_n=\O$, и $p(x,K_n)\ge 1-\varepsilon_n$}}\\
\text{{\it при $x\in K_{n+1}$, для всех $n\in N$.}}
\end{cases}
$\medskip

{\it Тогда для ЦМ существует инвариантная \underline{чисто конечно-аддитив}-\underline{ная} мера $\mu\in\Delta_{ba}$, причем $\mu(K_n)=1$ при всех $n\in N$.} \medskip

В [5] (Теорема 12.4) показывается, что это условие (Z) является также необходимым и достаточным для того, чтобы конечно-осредненная ЦМ (с любым параметром m), а значит и исходная ЦМ, \underline{не являлась} квазикомпактной.

\section{Примеры и методика их изучения}

Конечно-аддитивные меры в теории цепей Маркова возникают не только при общих фазовых пространствах. Они могут дать кое-что новое даже для ``почти'' феллеровских ЦМ на компакте. 

Напомним, что если $X$ топологическое пространство, $C(X)$ - пространство всех ограниченных непрерывных функций на $X$, и для оператора Т выполняется $T C(X)\subset C(X)$, то такая ЦМ называется феллеровской.  

Ниже мы рассматриваем пять взаимосвязанных простых ЦМ1--ЦМ5, заданных на отрезке $[0,1]$ с обычной борелевской сигма-алгеброй $\Sigma=\mathfrak{B}$ (сохраняем обозначение $\Sigma$), и даем некоторую \underline{методику} по использованию конечно-аддитивных мер при изучении ЦМ. \medskip

Схема возможных переходов для всех ЦМ одинаковая. Из любой начальной точки $x_0\in(0,1)$ возможны только два перехода: в точку $x_0^2$ с вероятностью $p(x_0)$, и в точку 0 (ноль) с вероятностью $1-p(x_0)$. Точки 0 (ноль) и 1 (единица) с вероятностью 1 переходят cами в себя, т.е. являются поглощающими (стационарными). ЦМ1 - ЦМ5 различаются только видом функции вероятности $p(x_0)$. Мы увидим, к какому непростому многообразию асимптотических (эргодических) свойств ЦМ приводят различные виды вероятностей $p(x_0)$. В работе автора [4] (§ 6, Пример 6.1) приводится только Пример 1 (ЦМ1), но без изучения его асимптотических и других свойств.

Приведем граф переходов (фазовый портрет) функционирования ЦМ1-ЦМ5 на Рисунке 1. У ЦМ3 и ЦМ5 ``буферы'' отсутствуют.

\begin{figure}[h!]
\center{\includegraphics[width=0.7\linewidth]{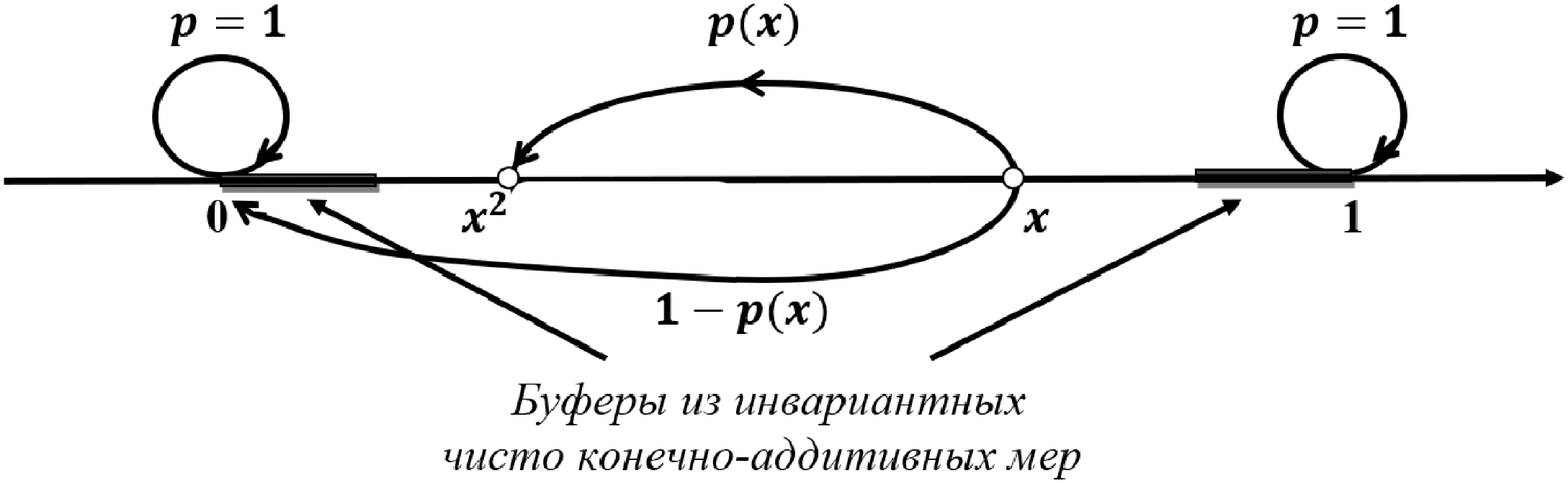}}
\caption{Фазовый портрет для ЦМ1 -- ЦМ5.}
\label{ris:image}
\end{figure}

Так как стационарная поглощающая точка $x=1$ недостижима из $[0,1)$ (``изолирована''), и нет циклических подклассов во всех пяти примерах, то, следуя Дубу ([8], Глава V, §5, случай l, стр. 196 русского издания), мы можем обойтись без использования средних по Чезаро (обычных в этой теории) от переходных вероятностей $p^n(x,E), n\in N$, т.е. без конечно-осредненных ЦМ, изучая лишь саму эту последовательность при $x\in[0,1)$. 
\medskip

{\bf\underline{Пример 1.}} Цепь Маркова, назовем ее ЦМ1, задается на отрезке $X=[0,1]$ следующими правилами.

Для любого $x\in (0,1)$ за один шаг возможны переходы - в точку $x^2$ с вероятностью $p(x)=x$, и в точку 0 (ноль) с вероятностью $1-p(x)=1-x$. Точки 0 (ноль) и 1 (единица) с вероятностью 1 переходят сами в себя, т.е.являются \underline{поглощающими} (стационарными). 

Формализуем переходную функцию такой ЦМ1. Для любого $x\in (0,1)$ определяем
$$p(x,\{x^2\})=x, p(x,\{0\})=1-x; p(0,\{0\})=1, p(1,\{1\})=1.$$

Все пять ЦМ мы будем изучать по одной и той же схеме, поэтому пронумеруем все дальнейшие шаги. \medskip

\begin{enumerate}
\item {\it ``Видимые'' инвариантные меры.} 

Легко видеть, что ЦМ1 имеет, ``по крайней мере'', \underline{две инва}-\underline{риантные сингулярные меры} - меры Дирака $\delta_0$ и $\delta_1$, сосредоточенные в точках 0 и 1 соответственно. Обе меры счетно-аддитивны. Другие возможные инвариантные меры ``невооруженным взглядом'' пока не просматриваются. \medskip

\item {\it Феллеровость.}

Пусть $C_{[0,1]}$ - банахово пространство всех непрерывных функций $f: [0,1]\to R$, и $f\in{C_{[0,1]}}$. Тогда для любого $x\in (0,1)$ выполняется 

$$g(x)=T{f(x)}=\int_{[0,1]}f(y)p(x,dy)=f(0)(1-x)+f(x^2)\cdot{x},$$
 
$$g(0)=T{f(0)}=f(0), g(1)=T{f(1)}=f(1).$$ 
Отсюда видно, что функция $g(x)=T{f(x)}$ непрерывна во всех точках $x\in [0,1]$, что означает, что ЦМ1 является \underline{фелле}\-\underline{ровской} на компакте [0,1].\medskip

\item {\it Асимптотика.}

Рассмотрим асимптотическое поведение ЦМ1. Ограничимся случаем, когда ЦМ1 ``запускается'' из произвольной фиксированной точки $x_0\in (0,1)$, которой соответствует начальная мера Дирака $\mu_0=\delta_{x_0}$, сосредоточенная в точке $x_0$. Начальная мера $\mu_0$ порождает марковскую последовательность счетно-аддитивных мер 
$$
\mu_n = A_{\mu_{n-1}}, n\in N, \mu_n(E)=p^n(x_0,E), E\in\Sigma, n\in N.
$$ 

Соответствующие им случайные величины $\xi_n$ могут принимать только два значения (два атома меры $\mu_n$): $\xi_n = x_0^{2^n}$ с вероятностью $\mu_n(\{x_0^{2^n}\})$, и $\xi_n=0$ с вероятностью $\mu_n(\{0\}), n\in N$. 

Прямым расчетом по индукции получаем следующие формулы:
$$
\mu_n(\{x_0^{2^n}\})=p^n(x_0,\{x_0^{2^n}\})=p(x_0, x_0^{2})\cdot p(x_0^{2}, x_0^{4}) \dots p(x_0^{2^{n-1}}, x_0^{2^n})=
$$ 
$$
=\prod_{k=0}^{n-1}x_0^{2^k}=x_0^{2^n-1}; \mu_n(\{0\})=1-x_0^{2^n-1}.
$$

Построим бесконечное ``траекторное множество'' $V_{x_0}=\{x_0$, $x_0^2, x_0^4,$ $x_0^8, \dots, x_0^{2^n}, \dots\}$ всех возможных ненулевых значений случайных величин $\xi_n$, порождаемых начальной точкой $x_0\in(0,1)$. 

Возьмем меру Дирака $\delta_0$, которая является инвариантной для ЦМ1: $A\delta_0=\delta_0$. Очевидно, $\delta_0(V_{x_0})=0$. Тогда при $n\in N$ выполняется 
$$
\mu_n(\{V_{x_0}\})-\delta_0(\{V_{x_0}\})|=\mu_n(\{V_{x_0}\})=\mu_n(\{x_0^{2^n}\})=x_0^{2^n-1},
$$
$$
\mu_n(\{0\})-\delta_0(\{0\})|=x_0^{2^n-1}.
$$
Поскольку других точек в ``носителях'' мер $\mu_n$ и $\delta_0$ нет, то для расстояния между ними получим: 
$$\rho(\mu_n,\delta_0)=\|\mu_n-\delta_0\|=\sup\limits_{E\in\Sigma}|\mu_n(E)-\delta_0(E)|=x_0^{2^n-1}\to0$$ для любого $x_0\in(0,1)$ при $n\to\infty$.

Итак, для любого начального $x_0\in(0,1)$ последовательность счетно-аддитивных марковских мер $\mu_n$ \underline{сходится} к инвариантной счетно-аддитивной мере $\delta_0$ \underline{в сильной метрической} \underline{топологии} в пространстве $ca(X,\Sigma_0)$, где $X=[0,1), \Sigma_0=\Sigma \cap [0,1)$ .

Однако, при любом $n\in N$ 
$$
\sup\limits_{0<x_0<1}\|\mu_n-\delta_0\|=\sup\limits_{0<x_0<1}x_0^{2^{n}-1}\equiv 1\nrightarrow 0 {\text{ при }} n\to\infty.
$$
Это означает, что \underline{сильная сходимость $\mu_n\to\delta_0$ не рав}\-\underline{номерна} по начальным точкам $x_0$. Тем более, эта сходимость не равномерна по всем возможным начальным мерам $\mu_0\in S_{ca}$ на $(0,1)$.

Отметим также, что если $x_0=0$, то $\mu_n\equiv\delta_0,n\in N$, т.е. формально $\|\mu_n-\delta_0\|\to0$ при $n\to\infty$; если $x_0=1$, то $\mu_n\equiv\delta_1, n\in N$ и $\|\mu_n-\delta_1\|\to 0$ при $n\to\infty$. \medskip

\item {\it Квазикомпактность.}

Если бы ЦМ1 была квазикомпактной, то сильная сходимость $\mu_n\to\delta_0$ была бы равномерной по $x_0\in(0,1)$. Следовательно, \underline{ЦМ1 не является квазикомпактной.}

Что этому мешает в Примере 1? Исходя из нашей Теоремы 1 и ее Следствия, можно предположить, что у ЦМ1 есть еще какие-то инвариантные меры, причем в большом количестве, которые мы не увидели при первом поверхностном рассмотрении в Пункте 1, т.е. ``неявные'' инвариантные меры. 

Изучением этой главной в статье задачи, применительно к ЦМ1, мы и займемся далее. \medskip

\item {\it Свойства ``неявных'' инвариантных мер, если они существуют.}

Пусть $\mu\in S_{ba}$ произвольная конечно-аддитивная вероятностная мера. Выполним следующие интегральные преобразования с оператором А цепи ЦМ1:
$$
A\mu(\{0\})=\int_{[0,1]}p(x,\{0\})\mu (dx)= \int_{\{0\}}+\int_{\{1\}}+\int_{(0,1)}=
$$
$$
=1\cdot\mu(\{0\})+0\cdot\mu(\{1\})+\int_{(0,1)}p(x,\{0\})\mu (dx)=
$$
$$
=\mu(\{0\})+\int_{(0,1)}(1-x)\mu (dx).
$$

Предположим теперь, что мера $\mu\in S_{ba}$ является инвариантной для ЦМ1, и она сингулярна с инвариантными мерами $\delta_0$ и $\delta_1$, т.е. $\mu((0,1))=1$. Тогда 
$$
\mu(\{0\})=A\mu(\{0\})=\mu(\{0\})+\int_{(0,1)}(1-x)\mu (dx),
$$
откуда получаем 
$$
\int_{(0,1)}(1-x)\mu (dx)=0.
$$

Это не ``равенство'', а ``уравнение'', которое мы должны решить относительно неизвестной меры $\mu$. На это уйдет несколько этапов. 

Пусть $0<a\le b<1$, т.е. $[a,b]\subset(0,1)$. Тогда 
$$
(1-b)\mu([a,b])\le \int_{[a,b]}(1-x)\mu(dx)\le(1-a)\mu([a,b])\le\mu([a,b]).
$$
Следовательно, если $\mu([a,b])>0$, то из левой части неравенств вытекает, что
$$
0=\int_{(0,1)}(1-x)\mu(dx)\ge\int_{[a,b]}(1-x)\mu(dx)>0,
$$
и мы получаем противоречие. Итак, для любого $[a,b]\subset(0,1)$, $\mu([a,b])=0$. Понятно, что для любого множества $E=\Sigma$, $E\subset[a,b]$ при произвольном $[a,b]\subset(0,1)$ выполняется $\mu(E)=0$. В частности, $\mu(\{a\})=0$ для любого $a\in(0,1)$.

Аналогично можно получить, что для любого $\varepsilon\in(0,1)$ верно $\mu((\varepsilon,1))=1$ и, соответственно, $\mu((0,\varepsilon))=0$. Это ключевое свойство таких мер.

Покажем теперь, что такие меры не счетно-аддитивны (хотя это и почти очевидно). Построим последовательность множеств 
$$
E_n=\left(1-\frac{1}{n},1\right), n\in N,
$$
для которой 
$$
E_1\supset E_2 \supset\dots, \lim E_n=\bigcap\limits_{n=1}^\infty E_n = \O.
$$
Тогда $\mu(E_n)=1, n\in N$, и 
$$
\lim\limits_{n\to\infty}\mu(E_n)=1\ne\mu(\lim E_n)=\mu(\O)=0.
$$ 
Это и означает, что мера $\mu$ не является непрерывной, т.е. не счетно-аддитивной. Итак, все конечно-аддитивные меры $\mu\in S_{ba}$, удовлетворяющие условиям $\mu((0,1))=1$ и $\int_{(0,1)}(1-x)\mu(dx)=0$, являются чисто конечно-аддитивными, и для них при любом $\varepsilon\in(0,1)$ верно $\mu((0,\varepsilon))=0$ и $\mu((\varepsilon,1))=1$, т.е. эти меры сосредоточены ``около единицы''. 

В общей теории конечно-аддитивных мер (см., например, [2]) доказывается, что такие меры существуют и их очень много, - мощность всего семейства таких мер ``около единицы'', - гиперконтинуум ($2^{2^{\aleph_0}}$).

Пока что мы всего лишь показали, что, \underline{если у ЦМ1 существу}-\underline{ет инвариантная} конечно-аддитивная \underline{мера} $\mu=A\mu$, $\mu((0,1))=1$, то она попадает в описанное выше семейство чисто конечно-аддитивных мер. \medskip

\item {\it Доказательство существования ``неявных'' инвариантных мер.}

Строим для ЦМ1 две последовательности $\varepsilon_n$ и $K_n$, описанные в условиях Теоремы 7.8.  [4], приведенной в пункте 2 настоящей статьи. Пусть $\varepsilon\in(0,1)$ произвольно и фиксировано, $\varepsilon_n=1-\varepsilon^{\frac{1}{2^n}}, n\in N$. Отсюда $\varepsilon^{\frac{1}{2^n}}=1-\varepsilon_n\to\varepsilon^0=1$ и $\varepsilon_n\to 0$ при $n\to\infty$.

Строим последовательность множеств $K_n\in\Sigma$: 
$$
K_1=(\varepsilon,1)=(\varepsilon^{\frac{1}{2^0}},1), K_2=(\sqrt\varepsilon,1)=(\varepsilon^{\frac{1}{2^1}},1),\dots,
$$
$$
K_n=(\varepsilon^{\frac{1}{2^{n-1}}},1), K_{n+1}=(\varepsilon^{\frac{1}{2^n}},1), \dots .
$$ 
Очевидно, $K_1\supset K_2\supset\dots$, $\lim K_n =\bigcap\limits_{n=1}^\infty K_n=\O$.

Пусть $x\in K_{n+1}$, т.е. $\varepsilon^{\frac{1}{2^n}}<x<1$. Тогда $p(x,K_n)=x>\varepsilon^{\frac{1}{2^n}}=1-\varepsilon_n$.

Теперь условия (критерий Z) Теоремы 7.8. [4] выполнены. Следовательно, для нашей ЦМ1 существует инвариантная чисто конечно-аддитивная мера $\mu\in\Delta_{ba}$, причем $\mu(K_n)=\mu((\varepsilon^{\frac{1}{2^{n-1}}},1))=1$ при $n\in N$. В частности $\mu((\varepsilon,1))=1$ и $\mu((0,\varepsilon))=0$ для любого $\varepsilon\in(0,1)$. Кроме того, $\mu(\{0\})=0$ и $\mu(\{1\})=0$. 

Итак, по ``крайней мере'', \underline{одна} инвариантная чисто конечно-аддитивная мера ``около единицы'' для ЦМ1 \underline{обнаружена}.
\medskip

\item {\it Количество ``неявных'' инвариантных мер.}

Покажем, что, на самом деле, таких инвариантных мер (линейно независимых) бесконечно много. Для нас психологически это важно, поскольку иначе наша Теорема 1 в Пункте 2 статьи окажется неверной.\medskip

Возьмем в качестве начальной точки нашей ЦМ1 трансцендентное число $x_0=1/\pi$ ($\pi\approx 3,14\dots$) и построим все порождаемое им ``траекторное множество'' $V_\pi$ точек вместе с прошлым и будущим, и добавим к нему точки 0 и 1. Счетное множество $V_\pi$ стохастически замкнуто, т.е. $p(x,V_\pi)=1$ для всех $x\in V_\pi$. Следовательно, мы можем сузить ЦМ1 с фазового пространства $([0,1],\Sigma)$ на фазовое пространство $(V_\pi,\Sigma_\pi)$, где $\Sigma_\pi=\Sigma\bigcap V_\pi=2^{V_\pi}$, с сохранением вида ее переходной функции. Получим новую частную ЦМ1. Все выявленные нами ранее свойства общей ЦМ1 (кроме феллеровости, которая здесь теряет смысл) сохраняются и для частной ЦМ1 на $V_\pi$. В частности, она имеет инвариантные счетно-аддитивные меры $\delta_0$ и $\delta_1$.

Не сложно (но громоздко) проверить, что для нее выполняются также условия (критерий Z) Теоремы 7.8. [4], влекущие существование для частной ЦМ1 инвариантной чисто конечно-аддитивной меры $\mu_\pi$ ``около единицы'', определенной только на счетном множестве $V_\pi$, и удовлетворяющей условию: $\mu_\pi(V_\pi\bigcap(\varepsilon,1))=1$ для всех $\varepsilon\in(0,1)$. Теорема 7.8. [4] и в этом случае применима, так как она сформулирована и доказана в [4] для любых измеримых пространств $(X,\Sigma)$ (``общая фаза'').

Теперь продолжим инвариантную меру $\mu_\pi$ с пространства $(V_\pi,\Sigma_\pi)$ на все исходное пространство $(X,\Sigma)$, доопределив ее ``нулями'' на всех множествах вне $V_\pi$. Несложно увидеть, что такая продолженная мера $\mu_\pi$ останется чисто конечно-аддитивной, и станет инвариантной для исходной общей ЦМ1 на $((0,1),\Sigma)$.\medskip
 
Теперь возьмем в качестве начальной точки нашей ЦМ1 другое трансцендентное число $1/e$ $(e\approx 2,71\dots)$ и построим его ``траекторное множество'' (вместе с точками 0 и 1) $V_e$. Очевидно, что для $V_e$ будут верны все рассуждения и выводы, сделанные выше для траекторного множества $V_\pi$. Его инвариантная чисто конечно-аддитивная мера $\mu_e$ также продолжается ``нулями'' до инвариантной меры на исходном пространстве для общей ЦМ1.
 
Числа $1/\pi$ и $1/e$ несоизмеримы, и порождаемые ими ``траекторные множества'' $V_\pi$ и $V_e$ не пересекаются. Следовательно, их инвариантные меры $\mu_\pi$ и $\mu_e$ сингулярны и линейно независимы, и их продолжения ``нулями'' на все пространство $([0,1],\Sigma)$ также сингулярны.\medskip
 
Итак, мы получили уже \underline{две} инвариантные чисто конечно-ад\-ди\-тив\-ные меры ЦМ1 ``около единицы'', сингулярные друг с другом. Множество всех трансцендентных попарно несоизмеримых чисел в интервале $(0,1)$ имеет мощность континуум $(2^{\aleph_0})$, и каждое такое число порождает свое траекторное множество со своей инвариантной чисто конечно-аддитивной мерой, попарно сингулярной и линейно независимой со всеми другими такими инвариантными мерами.
 
Мы показали, что наша ЦМ1 имеет \underline{не менее континуума} \underline{инвариантных} попарно сингулярных и линейно независимых чисто конечно-аддитивных мер ``около единицы''.\medskip

Выше, в Пункте 5, мы увидели, что ``около других мест'' в (0,1) инвариантных мер для ЦМ1 нет. В частности, нет инвариантных чисто конечно-аддитивных мер у ЦМ1 и ``около нуля'' с условием $\mu((0,\varepsilon))=1$ для всех $\varepsilon\in(0,1)$, хотя ЦМ1 и сходится к нулю. И этот факт весьма удивителен.\medskip

Вот мы и нашли причину, по которой мы не смогли воспользоваться нашей Теоремой 1 для ЦМ1. В условиях Теоремы требуется, чтобы $dim\Delta_{ba}=n<\infty$, а для ЦМ1 оказалось, что $dim\Delta_{ba}=\infty$, т.е. условия Теоремы 1 не выполнены. По этой же причине ЦМ1 и не обладает равномерной сходимостью и \underline{не является квазикомпактной}. 
 
\end{enumerate}  \medskip

{\bf\underline{Пример 2.}}
Берем за основу предыдущий Пример 1 с теми же возможными переходами на [0,1], и меняем местами вероятности этих переходов, т.е. $p(x)=1-x, 1-p(x)=x$. Полученную новую цепь назовем ЦМ2. Переходная функция ЦМ2 приобретает следующий вид: для любого $x\in(0,1)$ выполняется
$$
p(x,\{x^2\})=1-x, p(x,\{0\})=x; p(0,\{0\})=1, p(1,\{1\})=1.
$$

Проведем анализ ЦМ2 по такой же схеме, как и анализ ЦМ1, но сократив сопутствующие рассуждения.
\begin{enumerate} \medskip

\item {\it``Видимые'' инвариантные меры} у ЦМ2, такие же, как и у ЦМ1, -- две инвариантные счетно-аддитивные меры $\delta_0$ и $\delta_1$. Есть ли другие инвариантные меры, также пока не видно. \medskip    

\item {\it Феллеровость.} Пусть $f\in C_{[0,1]}$. Построив функцию $g(x)=T{f(x)}$, легко увидеть, что непрерывность $g(x)$ может нарушаться только в одной точке $x=1$, т.е. \underline{ЦМ2 не является} \underline{феллеровской} (назовем такие ЦМ \underline{``почти феллеровскими''}). \medskip

\item {\it Асимптотика.} Фиксируем произвольную начальную точку $x_0\in(0,1)$ с соответсвующей начальной мерой $\mu_0=\delta_{x_0}$. Строим марковскую последовательность мер $\mu_n=A\mu_{n-1},n\in N$, задающих случайные величины $\xi_n$, которые могут принимать лишь два значения с соответствующими вероятностями: $p\{\xi_n=x_0^{2^n}\}=\mu_n(\{x_0^{2^n}\})$,  $p\{\xi_n=0\}=\mu_n(\{0\}), n\in N$.

Прямым расчетом по индукции получаем следующие формулы:

$\mu_n(\{x_0^{2^{n}}\})=p^n(x_0,\{x_0^{2^{n}}\})=\prod\limits_{k=0}^{n-1} (1-x_0^{2^{k}});$

$\mu_n(\{0\})=1-\prod\limits_{k=0}^{n-1} (1-x_0^{2^{k}})$.

Так как ряд $\sum\limits_{k=0}^{\infty} x_0^{2^k}<\infty$, т.е. сходится, то и полученные произведения также монотонно \underline{убывая}, сходятся к \underline{положи-} \underline{тельному} числу 
$$
\lim\limits_{n\to\infty}\prod\limits_{k=0}^{n} (1-x_0^{2^k})=\prod\limits_{k=0}^{\infty} (1-x_0^{2^k})=\gamma(x_0), 0<\gamma(x_0)<1.
$$
Соответственно, и последовательность $\mu_n(\{0\})$ монотонно \underline{воз-} \underline{растает} и сходится к положительному числу $1-\gamma(x_0),$ $ 0<1-\gamma(x_0)<1$, т.е. $\mu_n(\{0\})$ не сходится к 0 и $\mu_n(\{0\})\ne\delta_0(\{0\})$ при любом начальном $x_0\in(0,1)$. Следовательно, расстояние $\rho(\mu_n,\delta_0)\nrightarrow 0$ при $n\to\infty$. \medskip

Итак, для любого начального $x_0\in(0,1)$ последовательность мер $\{\mu_n\}$ \underline{не сходится} к инвариантной мере $\delta_0$ \underline{в сильной} \underline{метрической топологии}. Более того, меры $\mu_n$ быстро и монотонно удаляются от меры $\delta_0$ до расстояния $\gamma(x_0)$. Однако, при этом положительные носители мер $\mu_n$, $x_0^{2^{n-1}}$, как числовая последовательность, быстро и монотонно сходятся к точке 0, т.е. к носителю меры $\delta_0$, в обычной евклидовой метрике числовой прямой. 

На всякий случай отметим, что $\rho(\mu_n,\delta_1)=\|\mu_n-\delta_1\|\equiv 1, n\in N$, т.е. последовательность $\mu_n$, естественно, не сходится сильно и ко второй инвариантной мере $\delta_1$.\medskip 

Если меры $\mu_n$ и $\delta_0$ рассматривать как элементы банахова пространства $ba([0,1],\Sigma)$, то последовательность $\mu_n$ \underline{не сходится} к $\delta_0$ также \underline{в слабой} и \underline{*-слабой топологии} поскольку существует множество $E=\{0\}$ такое, что $\mu_n(E)\nrightarrow\delta_0(E)$ при $n\to\infty$.

Если же рассматривать меры $\mu_n$ и $\delta_0$ как элементы банахова пространства $ca([0,1],\Sigma)$, то здесь нужно учесть, что пространство $ca([0,1],\Sigma)$ (равно $rca([0,1],\Sigma)$) является топологически сопряженным к пространству $C_{[0,1]}$, так как отрезок $[0,1]$ - компакт (см., например, [1]), поэтому, *-слабая топология в $ca([0,1],\Sigma)$ является хорошо известной $\Im_c$-топологией. И вот в этой $\Im_C$-тополо\-гии последовательность $\mu_n$ сходится (наконец-то!) к $\delta_0$. Несложно проверить, что при любом $x_0\in(0,1)$ и для всех $f\in C_{[0,1]}$ выполняется: 
$$
\langle f,\mu_n \rangle =\int_{[0,1]}f(x)\mu_n(dx)=
$$
$$
=f(x_0^{2^n})(1-x_0^{2^{n-1}})+f(0)\cdot x_0^{2^{n-1}}\to f(0)=\langle f,\delta_0\rangle, n\to\infty.
$$ 

Итак, для ЦМ2 при любом начальном $x_0\in(0,1)$ марковская последовательность мер $\mu_n$ \underline{сходится} к ивариантной мере $\delta_0$ только в \underline{$\Im_C$-слабой топологии.}

Понятны и тривиальные случаи при $x_0=0$ или $x_0=1$ (не будем комментировать). \medskip

\item {\it Квазикомпактность.} 

Из асимптотических свойств ЦМ2 следует, что она \underline{не является} \underline{квазикомпактной}. Следовательно, можно предположить, что, ЦМ2, так же, как и ЦМ1, имеет и другие, кроме $\delta_0$ и $\delta_1$, ``неявные'' инвариантные меры. \medskip
\item {\it Свойства ``неявных'' инвариантных мер.}

Проделаем интегральное преобразования, сделанные ранее для ЦМ1 в Пункте 5, с учетом перемены местами вероятностей переходов (опускаем). Получим следующее.

Если конечно-аддитивная мера $\mu$, с условием $\mu((0,1))=1$, является инвариантной для ЦМ2, то она удовлетворяет уравнению $\int_{(0,1)}x\cdot\mu(dx)=0$.

Отличие здесь от аналогичного уравнения для ЦМ1 в подынтегральной функции. Для ЦМ1 - это $f(x)=1-x$, для ЦМ2 - это $f(x)=x$. Проанализировав данное уравнение для ЦМ2 по той же схеме, что и для ЦМ1, получим следующее.

Если конечно-аддитивная мера $\mu, \mu((0,1))=1$, является решением интегрального уравнения для ЦМ2, то она обладает следующими свойствами: \medskip

1) $\mu([a,b])=0$ для любого $[a,b]\subset(0,1)$;

2) $\mu(\{\alpha\})=0$ для любого $a\in(0,1)$;

3) $\mu((0,\varepsilon))=1$ и $\mu((\varepsilon,1))=0$ для любого $\varepsilon\in(0,1)$.\medskip

Такие меры для ЦМ2 расположены ``около нуля'', а анологичные меры для ЦМ1 были расположены ``около единицы''.

Точно так же, как для ЦМ1, показываем, что любая такая мера является чисто конечно-аддитивной, т.е. не существует положительных счетно-аддитивных мер, удовлетворяющих соответствующему интегральному уравнению.

Если у ЦМ2 существует инвариантная конечно-аддитивная мера $\mu=A\mu$, $\mu((0,1))=1$, то она попадает в описанное выше семейство чисто конечно-аддитивных мер, расположенных ``около нуля''. \medskip      
 
\item {\it Доказательство существования ``неявных'' инвариантных мер.}

Доказательство для ЦМ2 проводится аналогично доказательству для ЦМ1, с использованием условия (критерия Z) Теоремы 7.8. [4], которое приведено выше в Пункте 2. Отличие лишь в том, что если для ЦМ1 мы строили последовательность множеств $K_n$ ``около единицы'', то теперь строим ее ``около нуля''.

Для произвольного $\varepsilon\in(0,1)$, берем 
$$
\varepsilon_n=\varepsilon^{\frac{1}{2^n}}, K_n=(0,\varepsilon^{2^{n-1}}), n\in N.
$$
Показываем, что $p(x,K_n)>1-\varepsilon_n$ для всех $x\in K_{n+1}$. Условия Теоремы 7.8. [4] выполнены.

Следовательно, для ЦМ2 существует инвариантная чисто конечно-аддитивная мера $\mu\in\Delta_{ba}$, причем $\mu((0,\varepsilon))=1$ для всех $\varepsilon\in(0,1)$. \medskip
\item {\it Количество ``неявных'' инвариантных мер.}

При произвольном $x_0\in(0,1)$, пораждаемые им траекторные множества для ЦМ1 и ЦМ2 совпадают. Поэтому для ЦМ2 можно повторить все рассуждения, приведенные в Пункте 7 для ЦМ1. Отличие только в том, что инвариантные чисто конечно-аддитивные меры для ЦМ1 удовлетворяют условию $\mu((\varepsilon,1))=1$ для $\varepsilon\in(0,1)$, а для ЦМ2 - условию $\mu((0,\varepsilon))=1$ для $\varepsilon\in(0,1)$. Но это никак не влияет на следующий общий вывод.

Наша ЦМ2 имеет не менее континуума инвариантных попарно  сингулярных и линейно-независимых чисто конечно-аддитивных мер, расположенных ``около нуля'', и других таких инвариантных мер, расположенных в других ``местах'' на $(0,1)$ - нет. Условия Теоремы 1 не выполнены.
\end{enumerate} 
\medskip

Для ЦМ2 справедливы замечания, которые мы сделали для ЦМ1 в конце текста из Примера 1.
\medskip

{\bf\underline{Пример 3.}}
Из двух ЦМ1 и ЦМ2 из Примеров 1 и 2 скомбинируем новую ЦМ3 с правилами переходов на $[0,1/2]$ такими же, как у ЦМ1 из Примера 1, а на $(1/2,1]$ - как у ЦМ2 из Примера 2. Запишем ее переходную функцию: 

Пусть $x\in[0,1/2]$. Тогда 
$$
p(x,\{x^2\})=x, p(x,\{0\})=1-x; p(0,\{0\})=1.
$$

Пусть $x\in(1/2,1]$. Тогда
$$
p(x,\{x^2\})=1-x, p(x,\{0\})=x; p(1,\{1\})=1.$$
Обе переходные функции совпадают при $x=1/2$, т.е. они непрерывно склеены в этой точке.

Такая комбинация переходных функций ЦМ1 и ЦМ2 сделана по следущим соображениям. У ЦМ1 не было инвариантных чисто конечно-аддитивных мер ``около нуля'', но были ``около единицы''. У ЦМ2, наоборот, - были такие меры ``около нуля'', но их не было ``около единицы''. В надежде на то, что мы получим ЦМ3 совсем без инвариантных чисто конечно-аддитивных мер и ``около нуля'', и ``около единицы'', мы и выбрали на левом промежутке $[0,1/2]$ переходную функцию от ЦМ1, а на правом промежутке $(1/2,1]$ - переходную функцию от ЦМ2 (это легко сделать, так как траекторные множества у ЦМ1 и ЦМ2 совпадают). 

Если наши надежды оправдаются, то ряд пунктов из нашей схемы изучения ЦМ1 и ЦМ2 станут излишними. Поэтому, для ЦМ3 мы сразу перейдем к главному. Вначале только отметим, что ЦМ3 имеет те же ``видимые'' инвариантные меры $\delta_0$ и $\delta_1$, и ее феллеровость нарушается лишь в одной точке $x=1$, т.е. она является почти феллеровской.   

Сделаем интегральные преобразования для ЦМ3 по аналогии с преобразованиями для ЦМ1 и ЦМ2. В результате, для любой инвариантной конечно-аддитивной меры $\mu$ ЦМ3 на $(0,1)$, т.е. при $\mu((0,1))=1$ получим следующее:

$$
0=\mu(\{0\})=A\mu(\{0\})=\int_{(0,1)}p(x,\{0\})\mu(dx)=
$$
$$
=\int_{(0,1/2]}(1-x)\mu(dx)+\int_{(1/2,1)}x\cdot\mu(dx).
$$

Поскольку оба интеграла справа неотрицательны, то они оба равны нулю: 
$$\int_{(0,1/2]}(1-x)\mu(dx)=0, \int_{(1/2,1)}x\cdot\mu(dx)=0.$$ 

Решаем первое уравнение.
$$
0=\int_{(0,1/2]}(1-x)\mu(dx)\ge\int_{(0,1/2]}\frac{1}{2}\cdot\mu(dx)=\frac{1}{2}\mu((0,1/2))\ge 0.
$$
Следовательно, $\mu((0,1/2])=0$. 

Решаем второе уравнение. 
$$
0=\int_{(1/2,1)}x\mu(dx)\ge\int_{(1/2,1)}\frac{1}{2}\cdot\mu(dx)=\frac{1}{2}\mu((1/2,1))\ge 0.
$$
Следовательно, $\mu((1/2,1))=0$, откуда вытекает, что $\mu((0,1))=0$. 

Отсюда следует, что ЦМ3 не имеет отличных от нулевых инвариантных конечно-аддитивных мер на интервале $(0,1)$, т.е. нет инвариантных мер ни счетно-аддитивных, ни чисто конечно-адди\-тивных, что и ожидалось.

Теперь, хотя это, возможно, и очевидно, формально проверим выполнение для ЦМ3 наших условий $(\alpha)$ и Теоремы 1 из настоящей статьи.

Выше мы доказали, что для ЦМ3 существует всего две инвариантные конечно-аддитивные меры $\delta_0$ и $\delta_1$, т.е. $\Delta_{ba}=co\{\delta_0,\delta_1\}$, и размерность, $dim\Delta_{ba}=2$, конечна. Меры $\delta_0$ и $\delta_1$ сингулярны, и их носители $K_1=\{0\}$ и $K_2=\{1\}$ стохастически замкнуты, таким образом, выполнены условия $(\alpha)$ и условия Теоремы 1. Следовательно, все инвариантные меры ЦМ3 счетно-аддитивны, т.е. $\Delta_{ba}\subset\Delta_{ca}$ (это только в нашей ЦМ3 это очевидно).

По Следствию из Теоремы 1, для ЦМ3 выполнены условия Дуба-Деблина (D), наше условие (*), и \underline{ЦМ3 является квазиком-\-} \underline{пактной}.\medskip

Подтвердим выполнение условий (D) (они приведены выше в начале Пункта 2). Возьмем меру $\varphi\in ca([0,1],\Sigma)$, $\varphi=(1/2)\delta_0+(1/2)\delta_1$, $\varepsilon=1/4, 1-\varepsilon=3/4, k=1$. 

Для множеств $E=\{0\}$ и $E=\{1\}$, $\varphi(E)=1/2>\varepsilon=1/4$, и рассматривать $p(x,E)$ не нужно. 

Пусть $E=(0,1)$. Тогда $\varphi((0,1))=0<\varepsilon=1/4$. Если $x\in(0,1/2]$, то $p(x,(0,1))=x\le 1/2<1-\varepsilon=3/4$. Если $x\in(1/2,1)$, то $p(x,(0,1))=1-x<1/2<1-\varepsilon$. Кроме того, $p(0,(0,1))=0<1-\varepsilon$, и $p(1,(0,1))=0<1-\varepsilon$. 

Для любого другого ``меньшего'' множества $E\subset(0,1)$, $E\in\Sigma$, тем более будет выполняться $\varphi(E)=0<\varepsilon$ и $p(x,E)<1-\varepsilon$ для всех $x\in[0,1]$. Итак, условие (D) действительно выполнено, и ЦМ3 является квазикомпактной.\medskip

Квазикомпактность ЦМ3 гарантирует сильную метрическую сходимость мер $\mu_n$ к инвариантным мерам $\delta_0$ и $\delta_1$ равномерно по начальным мерам $\mu_0$. Однако, покажем, как это выглядит на данном конкретном Примере 3.

Из оценок, полученных выше при проверке условия (D), для $\varepsilon=1/4$ при всех $x_0\in[0,1]$, $\mu_0=\delta_{x_0}$, выполняется 
$$
\mu_1((0,1))=A\mu_0((0,1))=p^1(x_0,(0,1))<1-\varepsilon.
$$ 
Далее, при $x_0\in(0,1)$, делаем следующие преобразования и оценки: 
$$
\mu_2((0,1))=A\mu_1((0,1))=p^2(x_0,(0,1))=
$$
$$
=\int_{[0,1]}p^1(y,(0,1)) p^1(x_0,dy)=\int_{\{0\}}+\int_{\{1\}}+\int_{(0,1)}=
$$
$$
=0+0+\int_{(0,1)}p^1(y,(0,1))p^1(x_0,dy)<(1-\varepsilon)p^1(x_0,(0,1))<(1-\varepsilon)^2.
$$

Тривиальные случаи для $x_0=0$ и $x_0=1$ (их следует рассматривать отдельно) опускаем. 

По индукции получаем, что $\mu_n((0,1))<(1-\varepsilon)^n$ при $n\in N$ и для всех $x_0\in[0,1]$. Следовательно, $\lim\mu_n((0,1))=0$ при $n\to\infty$ равномерно по $x_0\in[0,1]$, и с экспоненциальной скоростью сходимости. 

Отсюда выводим, что при $x_0\in[0,1), x_0\ne 1$, выполняется $\mu_n(\{0\})\ge 1-(1-\varepsilon)^n$. Следовательно, $\lim\mu_n(\{0\})=1$ при $n\to\infty$, равномерно по $x_0\in[0,1)$, и с экспоненциальной скоростью.

При $x_0=1$ имеем $\mu_n(\{1\})\equiv\delta_1(1)=1, n\in N$, т.е., формально, $\mu_n(\{1\})$ ``сходится'' к 1, ``равномерно'' по единственной начальной точке $x_0=1$, из которой можно попасть в носитель меры $\delta_1$, т.е. в точку $x=1$.

Опускаем простые детали и получаем главный вывод. 

Марковская последовательность мер $\mu_n$ для ЦМ3 при $x_0\in[0,1)$ сходится в сильной метрической топологии к инвариантной мере $\delta_0$, т.е. $\rho(\mu_n,\delta_0)=\|\mu_n-\delta_0\|\to 0$, равномерно по всем начальным точкам $x_0\in[0,1)$ и экспоненциально быстро. Можно показать, что эта сходимость равномерна и по всем начальным вероятностным счетно-аддитивным мерам $\mu_0$, а не только по начальным мерам Дирака $\mu_0=\delta_{x_0}$.

При $x_0=1, \rho(\mu_n,\delta_1)\to 0$ также ``равномерно'' по $x_0=1$.

Различие в свойствах инвариантных мер $\delta_0$ и $\delta_1$ можно охарактеризовать следующим образом: мера $\delta_0$ устойчива (притягивающая), а мера $\delta_1$ неустойчива (отталкивающая, изолированная).

В Примерах 1 и 2 сильной равномерной сходимости к $\delta_0$ ЦМ1 и ЦМ2 препятствовали мощные \underline{буферы} из инвариантных чисто конечно-адди\-тив\-ных мер, прилипших в $\Im_C$-топологии к инвариантным счетно-адди\-тив\-ным мерам $\delta_1$ и $\delta_0$, соответственно.\medskip

{\bf\underline{Пример 4.}} Также, как и в Примере 3, скомбинируем из двух ЦМ 1 и ЦМ2 новую ЦМ4, но в другом порядке. На отрезке $[0,1/2]$ возьмем правило перехода из ЦМ2, а на полуинтервале $(1/2,1]$ берем правило перехода из ЦМ1. Не будем расписывать очевидные формулы для новой переходной функции у ЦМ4.

У предыдущей ЦМ3 мы взяли на указанных промежутках ``лучшие'' для сходимости переходные функции из ЦМ1 и ЦМ2. Теперь же, для ЦМ4 мы берем на этих промежутках ``худшие'' переходные функции из ЦМ2 и ЦМ1, ``тормозящие'' сходимость. Интуитивно ясно, что ЦМ4 унаследует ``худшие'' для сходимости свойства ЦМ2 и ЦМ1 (на наш взгляд, эти ``худшие'' свойства сложнее и более интересны).

Предлагаем читателю самому исследовать свойства ЦМ4 по схеме изучения ЦМ1 и ЦМ2 и получить следующее. ЦМ4 обладает двумя мощными \underline{буферами} из инвариантных чисто конечно-аддитивных мер около нуля и около единицы, имеет две инвариантные счетно-аддитивные меры, и не является квазикомпактной. Однако, при любом начальном $x_0\in(0,1)$ марковская последовательность счетно-аддитивных мер $\mu_n=A\mu_{n-1}$  $\Im_C$-слабо сходится к инвариантной счетно-аддитивной мере Дирака $\delta_0$, но не сходится к $\delta_0$ в метрической топологии. \medskip

{\bf\underline{Пример 5.}} Рассмотрим туже общую схему задания цепи Маркова на отрезке $[0,1]$, что и в ЦМ1- ЦМ4, т.е. с возможными переходами для любой $x_0\in[0,1)$ в точки $x_0^2$ и 0, но максимально упрощаем вероятности этих переходов. Полагаем, что $p(x_0,\{x_0^2\})=p$ и $p(x_0,\{0\})=q$, где вероятности $p$ и $q$ не зависят от $x_0\in[0,1)$ и $p+q=1$. По-прежнему полагаем $p(0,\{0\})=1$ и $p(1,\{1\})=1$. Построенную новую цепь Маркова обозначим ЦМ5. Все наши ЦМ1 - ЦМ5 однородны по времени, а новая ЦМ5 еще однородна и по пространству $[0,1)$ в указанном смысле.

Несмотря на максимальные упрощения, ЦМ5 не такая уж и ``хорошая'' - она всего лишь ``почти феллеровская''. Непрерывность функций $g(x)=Tf(x)$ для непрерывных $f(x)$ может нарушаться в точке $x=1$ (при $q\ne 0$). \medskip

Применяя нашу схему изучения ЦМ1 и ЦМ2 несложно показать следующее.

При любом $p, q, 0<p,q<1,$ у ЦМ5 есть всего две инвариантные счетно-аддитивные меры Дирака $\delta_0$ и $\delta_1$ и отсутствуют другие инвариантные конечно-аддитивные меры, т.е. $dim\Delta_{ba}=2<\infty$. Выполнены условия нашей Теоремы 1 со всеми вытекающими следствиями. 

Покажем лишь выполнение условий Дуба-Деблина (D). Берем меру $\varphi\in ca(X,\Sigma)$ как "средневзвешенную" из двух инвариантных мер Дирака $\delta_0$ и $\delta_1$: $\varphi=p\delta_0+q\delta_1$. Полагаем $\varepsilon=(1/2)\min(p,q)$ и $k=1$. Тогда для множеств $E=\{0\},  E=\{1\}$ выполняется $\varphi(\{0\})=p>\varepsilon$, и $\varphi(\{1\})=q>\varepsilon$, следовательно, эти множества мы снимаем с дальнейшего рассмотрения. Для множества $E=(0,1)$ выполняется $\varphi((0,1))=0<\varepsilon$. Оцениваем переходную функцию на этом множестве. Для трех возможных вариантов получаем: $p(0,(0,1))=0<1-\varepsilon$, $p(1,(0,1))=0<1-\varepsilon$, и при всех $x\in (0,1)$ выполняется $p(x, (0,1))=p<1-\varepsilon$. Для всех других множеств $E \subset (0,1)$, $E\in\Sigma$, тем более выполняется $p(x,E)<1-\varepsilon$ для всех $x\in[0,1]$. Следовательно, условия квазикомпактности (D) для ЦМ5 выполнено.

Эргодические свойства упрощенной ЦМ5 легко просматриваются и без использования чисто конечно-аддитивных мер. Но мы приводим этот Пример 5 для того, чтобы показать, что вся эта простота возможна лишь потому, что у ЦМ5 отсутствуют инвариантные чисто конечно-аддитивные меры. И они могут появляться, как только мы начнем варьировать параметры $p,q$ у ЦМ5 в зависимости от точки (состояния) $x$, что приведет к резкому изменению асимптотических и качественных свойств ЦМ.\medskip 

Мы сознательно выбрали простые Примеры 1-5, в которых инвариантные меры $\delta_0$ и $\delta_1$ для ЦМ 1-5 имеют одноточечные носители, что исключает существование циклических эргодических подклассов у квазикомпактных ЦМ. Свою важную функцию конечно-аддитивные меры, которым и посвящена настоящая статья, выполняют не внутри, а снаружи носителей инвариантных счетно-аддитивных мер, если последние существуют.   
\medskip

Сделаем общий неформальный комментарий к проведенному анализу ЦМ 1-5. Инвариантные счетно-аддитивные меры для ЦМ, сами по себе, дают очень мало информации об асимптотических свойствах ЦМ. Они могут вообще отсутствовать, но ЦМ при этом может обладать хорошей асимптотикой, и сходиться в каком-либо смысле к неинвариантной мере, или к чисто конечно-аддитивной мере (такая ситуация изучается в работе автора [5], Теоремы 13.1 и 13.2). В Примерах 1-5 показано, что большую ``ответственность'' за асимптотические свойства ЦМ несут инвариантные чисто конечно-аддитивные меры. А только один факт их отсутствия уже обеспечивает для ЦМ максимально хорошую эргодичность. Но это, скорее ``вырожденный случай'', чем ``норма''.\medskip   

Заметим также, что если более внимательно посмотреть на ``плохой'' Пример 2, то он окажется намного интереснее, чем ``хорошие'' Примеры 3, 4, 5  (но это уже выходит за рамки проблемы, изучаемой в теоретической части статьи, напишем об этом подробнее позже).\medskip

И еще одно общее замечание. При работе с чисто конечно-аддитив\-ны\-ми мерами вызывает затруднения то обстоятельство, что они являются продуктом ``трансфинитной индукции'', т.е. неконструктивной Аксиомы выбора. Существование таких мер из Аксиомы выбора следует, а вот наличие исчерпывающего и \underline{внятного} \underline{описания одной конкретной} чисто конечно-аддитивной меры не гарантируется, для них даже отсутствует понятие ``носителя'' меры. Но можно локализовать с той или иной степенью точности ``местоположение'' чисто конечно-аддитивной меры, что делает их все более и более понятными. Именно такой подход и иллюстрируют наши Примеры 1, 2 и 4. Но, не следует забывать, что чисто конечно-аддитивные меры могут быть устроены намного сложнее, чем те, которых породили ЦМ1, ЦМ2 и ЦМ4 в Примерах 1, 2 и 4.   

\section{Заключение}

В работе представлен новый результат по теории общих цепей Маркова, продолженных на пространство конечно-аддитивных мер. Полученный факт открывает новые возможности по дальнейшим исследованиям в данном направлении. Рассмотрены примеры на отрезке [0,1] с демонстрацией техники работы с чисто конечно-аддитивными мерами. 

Основные положения настоящей работы были предварительно кратко анонсированы (без доказательств) в [12], и опираются на теоремы работ автора [4], [5].

\end{document}